\documentstyle[12pt]{article}
\input xy
\xyoption{all}

\makeindex

\newcommand{\cof}{\mathrm{cof\kern1.2pt}}

\newcommand{\vsp}{\vskip 1em}

\newcommand{\rank}{{\mbox{rank\kern1.3pt}}}

\newcommand{\trace}{{\mbox{trace\kern1.3pt}}}

\newcommand{\bt}{\begin{theorem}}
\newcommand{\et}{\end{theorem}}
\newtheorem{lemma}{Lemma}
\newtheorem{theorem}[lemma]{Theorem}
\newtheorem{corollary}[lemma]{Corollary}

\newcommand{\del}{\delta}
\newcommand{\tms}{\times}

\newcommand {\be}{\begin{equation}}
\newcommand {\ee}{\end{equation}}

\def \qed {\hfill \vrule height6pt width 6pt depth 0pt}
\def \qed {\hfill \vrule height6pt width6pt depth0pt}

\def \textit{\it}
\def \bt{\begin{theorem}}
\def \et{\end{theorem}}
\def \bl{\begin{lemma}}
\def \el{\end{lemma}}
\def \bc{\begin{corollary}}
\def \ec{\end{corollary}}
\def \be{\begin{equation}}
\def \ee{\end{equation}}

\def \text{\mbox}

\title{Squared distance matrix of a weighted tree}
\author{Ravindra B. Bapat \\ Indian Statistical
Institute \\ New Delhi, 110016, India \\e-mail:
rbb@isid.ac.in}
\begin{document}
\maketitle

\begin{abstract}
Let $T$ be a tree with vertex set $\{1, \ldots, n\}$ such that each edge is assigned a nonzero
weight. The squared distance matrix of $T,$
denoted by $\Delta,$ is the $n \times n$ matrix with $(i,j)$-element $d(i,j)^2,$
where $d(i,j)$ is the sum of the weights of the edges on the $(ij)$-path. We obtain 
a formula for the determinant of $\Delta.$ 
A formula for $\Delta^{-1}$ is also obtained, under certain conditions. The results
generalize known formulas for the unweighted case.

\vspace{3mm}
\end{abstract}

\section{Introduction}

Let $G$ be a connected graph with vertex set $V(G) = \{1, \ldots, n\}.$
The distance between vertices $i,j \in V(G)$, denoted $d(i,j),$ is the minimum length
(the number of edges) of a path from $i$ to $j$ (or an $ij$-path). We set $d(i,i) = 0,
i = 1, \ldots, n.$ The distance matrix $D(G),$ or simply $D,$ is the $n \times n$
matrix with $(i,j)$-element $d_{ij} = d(i,j).$

A classical result of Graham and Pollak \cite{pollak} asserts that if $T$ is a tree
with $n$ vertices, then the determinant of the distance matrix $D$ of $T$ is
$(-1)^{n-1}(n-1)2^{n-2}.$ Thus the determinant depends only on the number of vertices 
in the tree and not on the tree itself. A formula for the inverse of the distance matrix of a tree
was given by Graham and Lov\'asz  \cite{lovasz}. Several extensions and generalizations 
of these results have been proved (see, for example \cite{bapat},
\cite{kirkneum}, \cite{lalpati}, \cite{yeh}, \cite{ding} and the references contained therein).

Let $T$ be a tree with vertex set $\{1, \ldots, n\}$ and let $D$ be the distance matrix
of $T.$ The squared distance matrix $\Delta$ is defined to be the Hadamard product
$D \circ D,$ and thus has the $(i,j)$-element $d(i,j)^2.$ A formula for the determinant of $\Delta$ was proved in \cite{krish1}, while the inverse and the inertia of $\Delta$ were considered in \cite{krish2}.

In this paper we consider weighted trees.
Let $T$ be a tree with vertex set $V(T) = \{1,\ldots,n\}$ and edge set
$E(T) = \{e_1, \ldots, e_{n-1}\}.$ We assume that each edge is assigned a
weight and let the weight assigned to $e_i$ be denoted $w_i,$ which is
a nonzero real number (not necessarily positive).

For $i,j \in V(T), i \not = j,$ the distance $d(i,j)$ is defined to be the sum
of the weights of the edges on the (unique) $ij$-path. We set $d(i,i)
= 0, i=1, \ldots,n.$ Let $D$
be the $n \tms n$ distance matrix with $d_{ij} = d(i,j).$

The Laplacian of $T$ is the $n \tms n$ matrix defined as follows. The rows and
the columns of $L$ are indexed by $V(T).$ For $i \not = j,$ the $(i,j)$-element
is $0$ if  $i$ and $j$ are not adjacent. If $i$ and $j$ are adjacent,
and if the edge joining them is $e_k,$ then the $(i,j)$-element of $L$
is set equal to $-1/w_k.$ The diagonal elements of $L$ are defined so that
$L$ has zero row (and column) sums.

The paper is organized as follows.
In this section we review some basic properties of the distance matrix of a tree
such as formulas for its determinant and inverse. Some preliminary
results are obtained in Section 2. Sections 3 and 4 are devoted to the determinant
and the inverse of $\Delta,$ respectively.

\vsp

\noindent 
{\bf Example} Consider the tree
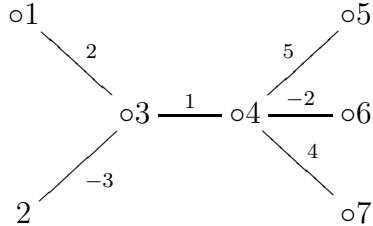
\[\xymatrix{\circ \mbox{1} \ar@{-}[rd]^{2} & & & \circ \mbox{5}  & & \\
& \circ \mbox{3}\ar@{-}[ld]^{-3}\ar@{-}[r]^{1} &\circ \mbox{4} \ar@{-}[ru]^{5} \ar@{-}[r]^{-2} \ar@{-}[rd]^{4}& \circ \mbox{6}  & \\
\mbox{2}& & & \circ \mbox{7}  & &
}\]

The Laplacian of the tree is given by
$$\left [ \begin{array}{rrrrrrr}
1/2 & 0 & -1/2 & 0 & 0 & 0 & 0\\
0 & -1/3 & 1/3 & 0 & 0 & 0 & 0\\
-1/2 & 1/3 & 7/6 & -1 & 0 & 0 & 0\\
0 & 0 & -1 & 19/20 & -1/5 & 1/2 & -1/4\\
0 & 0 & 0 & -1/5 & 1/5 & 0 & 0\\
0 & 0 & 0 & 1/2 & 0 & -1/2 & 0\\
0 & 0 & 0 & -1/4 & 0 & 0 & 1/4\\
\end{array} \right ].$$

\vskip 1em

We let $Q$ be the $n \tms (n-1)$ vertex-edge incidence matrix of the underlying
unweighted tree, with an orientation assigned to each edge.
Thus the rows and the columns of $Q$ are indexed by $V(T)$ and $E(T)$ respectively.
If $i \in V(T), e_j \in E(T),$ the $(i,j)$-element
of $Q$ is $0$ if $i$ and $e_j$ are not incident, it is $1 (-1)$
if $i$ and $e_j$ are incident and $i$ is the initial (terminal) vertex of $e_j.$
It is well-known \cite{bapat} that $Q$ has rank $n-1$ and  any minor of $Q$ is either
$0$ or $\pm 1$ (thus $Q$ is totally unimodular).

Let $F$ be the $n \tms n$ diagonal matrix with diagonal elements $w_1, \ldots, w_{n-1}.$
It can be verified that $L = QF^{-1}Q'.$

\begin{lemma} \label{lem31}
The following assertions are true:
\begin{description}
\item {(i)} $Q'DQ = -2F.$

\item {(ii)} $LDL = -2L.$

\end{description}
\end{lemma}

\noindent  {\bf Proof}
(i).  The result follows from the following observation which
is easily verified: If $e_p = \{i,j\}$ and $e_q=\{k, \ell\}$
are edges of $T,$ then
$$d(i,k) + d(j, \ell) - d(i,\ell) - d(j,k) $$
equals $0$ if $e_p$ and $e_q$ are distinct, and equals
$- 2w_p,$ if $e_p = e_q.$

(ii). We have
\begin{eqnarray*}
LDL &=& QF^{-1}Q'DQF^{-1}Q'\\
& = & QF^{-1}(-2F)F^{-1}Q' {\mbox{  by (i) }}\\
&=& -2 QF^{-1}Q' \\
&=& -2L,
\end{eqnarray*}
and the proof is complete. \qed

\vskip 1em

Let $\delta_i$ denote the degree of the vertex $i, i = 1, \ldots, n,$
and let $\delta$ be the $n \tms 1$ vector with components
$\delta_1, \ldots, \delta_n.$ We set $\tau_i
= 2 - \delta_i, i = 1, \ldots, n,$ and let $\tau$ be the $n \tms 1$
vector with components $\tau_1, \ldots, \tau_n.$

\begin{theorem}\label{lem32}
The following assertions are true:
\begin{description}
\item {(i)} $\det D = (-1)^{n-1}2^{n-2}(\sum_{i} w_i)(\prod_i w_i).$

\item {(ii)} If $\sum_i w_i \not = 0,$ then $D$ is nonsingular and
$$D^{-1} = -\frac{1}{2}L + \frac{1}{2\sum_i w_i}\tau \tau'.$$

\item {(iii)} $D \tau = (\sum_i w_i){\bf 1}.$
\end{description}

\end{theorem}

\noindent  {\bf Proof} Parts (i) and (ii) are well-known, see for example,
\cite{kirkneum}. To prove (iii), note that from (ii),
$$D^{-1} {\bf 1} = \frac{1}{2\sum_i w_i} \tau \tau' {\bf 1}
= \frac{1}{\sum_i w_i} \tau,$$
since ${\bf 1}'\tau = 2.$ It follows that $D \tau = (\sum_i w_i){\bf 1}$
and the proof is complete. \qed

\section{Preliminary results}

We now turn to the main results for the case of a weighted tree. 
Let $T$ be a tree with vertex set $V(T) = \{1,\ldots,n\}$ and edge set
$E(T) = \{e_1, \ldots, e_{n-1}\}.$ Let $w_1, \ldots, w_{n-1}$ be the
edge-weights.
Recall that $\delta_i$ is the degree of vertex $i$ and $\tau_i = 2 - \delta_i.$
We write $j \sim i$ if vertex $j$ is adjacent to vertex $i.$
We let $\hat{\delta}_i$ be the weighted degree of $i,$ which is defined as
$$\hat{\delta}_i = \sum_{j:j \sim i} w(\{i,j\}), i = 1, \ldots, n.$$ Let $\hat{\delta}$
be the $n \tms 1$ vector with components $\hat{\delta}_1, \ldots, \hat{\delta}_n.$

Let $\Delta$ be the squared distance matrix of $T,$ which is the $n \tms n$
matrix with its $(i,j)$-element equal to $d_{ij}^2$ or equivalently, $d(i,j)^2.$
The next result was obtained in \cite{krish2} for the unweighted case,

\begin{lemma} \label{lem33}
$\Delta \tau = D \hat{\delta}.$
\end{lemma} 

\noindent  {\bf Proof}
Let $i \in \{1, \ldots, n\}$ be fixed. For $j \not = i,$
let $\gamma(j)$ be the predecessor of $j$ on the $ij$-path (in the underlying unoriented tree). Let $e^j$ be the edge 
$\{\gamma(j),j\}$ and set $\theta^j = \hat{\delta}_j - w(e^j).$
We have
\begin{eqnarray}
& & 2 \sum_{j=1}^n d(i,j)^2 \nonumber \\
&=& \sum_{j=1}^n d(i,j)^2 + \sum_{j \not = i}
(d(i,\gamma(j)) + w(e^j))^2 \nonumber \\
&=& \sum_{j=1}^n d(i,j)^2 + \sum_{j \not =i} d(i,\gamma(j))^2 
+ 2 \sum_{j \not =i} d(i,\gamma(j))w(e^j) +  \sum_{j \not =i} w(e^j)^2. \label{eq31}
\end{eqnarray}

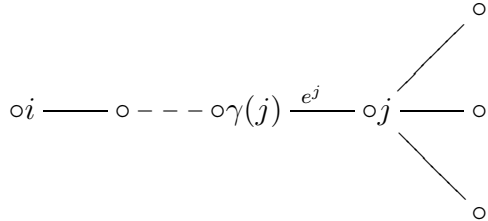
\[\xymatrix{
  &  &  & & \circ\\
\circ \mbox{$i$} \ar@{-}[r] & \circ \ar@{--}[r] & \circ \ar@{-}[r]^{e^j} \mbox{$\gamma(j)$} & 
\circ \mbox{$j$} \ar@{-}[r] \ar@{-}[ru] \ar@{-}[rd]
& \circ \\
 &  &  & & \circ
}\]

Note that
\begin{equation} \label{eq32}
\sum_{j \not =i} d(i,\gamma(j))^2  = \sum_{j=1}^n (\delta_j-1)d(i,j)^2,
\end{equation}
since vertex $j$ serves as a predecessor of $\delta_j - 1$ vertices in
paths from $i.$
Also note that
\begin{equation} \label{eq33}
\sum_{j \not =i} w(e^j)^2  = \sum_{k=1}^{n-1} w(e_k)^2.
\end{equation}
We have
\begin{eqnarray}
& &  \sum_{j=1}^n d(i,j)\hat{\delta}_j \nonumber \\
&=& \sum_{j \not = i} (d(i,\gamma(j) + w(e^j))(w(e^j) + \theta^j) \nonumber \\
&=& \sum_{j \not = i} d(i,\gamma(j))w(e^j)  + \sum_{j \not =i} w(e^j)^2 
+  \sum_{j \not =i} (d(i,\gamma(j)) + w(e^j))\theta^j. \label{eq34}
\end{eqnarray} 
Observe that  $\theta^j$ is the sum of the weights of all the edges incident to $j,$
except the edge $e^j,$ which is on the $ij$-path. Thus  $(d(i,\gamma(j)) + w(e^j))\theta^j$
equals $\sum d(i, \gamma(\ell))w(e^\ell),$ where the summation is over all vertices
adjacent to $j,$ except $i.$ Therefore it follows that
\begin{equation} \label{eq35}
\sum_{j \not = i} d(i,\gamma(j))w(e^j)  =  \sum_{j \not =i} (d(i,\gamma(j)) + w(e^j))\theta^j. 
\end{equation}
From (\ref{eq31})-(\ref{eq35}) we get
$$2 \sum_{i=1}^n d(i,j)^2 = \sum_{j=1}^n d(i,j)^2 \delta_j + \sum_{j=1}^n d(i,j)\hat{\delta}_j,$$
which is equivalent to
$$ \sum_{i=1}^n d(i,j)^2 \tau_j =  \sum_{j=1}^n d(i,j)\hat{\delta}_j,$$
and the proof is complete. \qed
 
\vsp

Next we define the edge orientation matrix of $T.$ We assign an orientation to each 
edge of $T.$ Let $e_i = (p, q); e_j = (r, s)$ be edges of $T.$ We say
that $e_i$ and $e_j$  are similarly oriented, denoted by $e_i \Rightarrow e_j,$   if 
$d(p, r) = d(q, s).$  Otherwise $e_i$ and $e_j$  are
said to be oppositely oriented, denoted by $e_i \rightleftharpoons e_j.$ 
 For example,  in the following diagram $e_i$ and $e_j$  are
similarly oriented.
\[\xymatrix{
\circ \mbox{$p$} \ar[r] & \circ \mbox{$q$} \ar@{--}[r] & \circ \ar[r] \mbox{$r$} & \circ \mbox{$s$}  
}\]

The edge orientation matrix of T is the $(n-1) \times (n-1)$ matrix $H$ having the rows and the columns
indexed by the edges of $T.$ The $(i, j)$-element of $H,$  denoted by $h(i, j)$
 is defined to be $1(-1)$ if the corresponding edges $e_i,e_j$ of  $T$ are similarly (oppositely) oriented.
The  diagonal elements  of $H$ are set to be $1.$ We assume that the same orientation is used
while defining the matrix $H$ and the incidence matrix $Q.$

If the tree $T$ has no vertex of degree $2,$ then we let $\hat{\tau}$ 
be the diagonal matrix with diagonal elements $1/\tau_1, \ldots,  1/\tau_n.$
We state some basic properties of $H$ next, see \cite{krish1}.

\begin{theorem}  \label{lem331}
Let $T$ be a directed tree on $n$ vertices, let $H$ and $Q$ be
the edge orientation matrix and the vertex-edge incidence matrix of $T,$ respectively.
Then  $\det H = 2^{n-2} \prod_{i=1}^n \tau_i.$ Furthermore, if $T$ has no vertex
of degree $2,$ then $H$ is nonsingular and 
$H^{-1} = \frac{1}{2} Q'\hat{\tau}Q.$
\end{theorem}

Let $w_1, \ldots, w_{n-1}$ be the edge-weights. Recall that $F$ be the diagonal matrix
with diagonal elements $w_1, \ldots, w_{n-1}.$

Also note that,
$$(FHF)_{ij}  =
\left \{ \begin{array}{ll}
w_iw_j & \mbox{ if } e_i \Rightarrow e_j\\
-w_iw_j &  \mbox{ if } e_i \rightleftharpoons  e_j.
\end{array} \right.
$$

\begin{lemma} \label{lem34}
$Q' \Delta Q = -2 FHF.$
\end{lemma}

\noindent  {\bf Proof} 
For $i,j \in \{1, \ldots, n-1\},$ let the edge $e_i$ be from $p$ to $q$
and the edge $e_j$ be from $r$ to $s.$ Then
\be \label{eq199}
(Q'\Delta Q)_{ij} =
\left \{ \begin{array}{ll}
d(p,r)^2 + d(q,s)^2 - d(p,s)^2 - d(q,r)^2 & \mbox{ if } e_i \Rightarrow e_j \\
d(p,s)^2 + d(q,r)^2 - d(p,r)^2 - d(q,s)^2 & \mbox{ if } e_i \rightleftharpoons e_j
\end{array} \right.
\ee
Let $d(r,s) = \alpha.$ It follows from (\ref{eq199}) that
\begin{eqnarray*} 
(Q'\Delta Q)_{ij} &=&
\left \{ \begin{array}{ll}
(w_i+\alpha)^2 + (w_j+\alpha)^2 - (w_i+w_j+\alpha)^2 - \alpha^2 = -2w_iw_j & \mbox{ if } e_i \Rightarrow e_j \\
(w_i+w_j+\alpha)^2 + \alpha^2 - (w_i+\alpha)^2 - (w_j+\alpha)^2 = 2w_iw_j & \mbox{ if } e_i \rightleftharpoons e_j
\end{array} \right. \\
&=& -2(FHF)_{ij},
\end{eqnarray*}
and the proof is complete.   \qed

\vsp

Let $\tilde{\tau}$ be the diagonal matrix with diagonal elements
$\tau_1, \ldots, \tau_n.$

\begin{lemma} \label{lem35}
$ \Delta L = 2 D \tilde{\tau} - {\bf 1} \hat{\delta}'.$
\end{lemma}

\noindent  {\bf Proof}
Let $i,j \in \{1, \ldots, n\}$ be fixed. Let vertex $j$ have degree $p.$ Suppose $j$
is adjacent to vertices $u_1, \ldots, u_p$ and let $e_{\ell_1}, \ldots, e_{\ell_p}$
be the corresponding edges with weights $w_{\ell_1}, \ldots, w_{\ell_p},$
respectively. We consider two cases.

\noindent
{\it Case (i).} $i = j.$ We have
\begin{eqnarray*}
(\Delta L)_{jj} &=& \sum_{k=1}^n d(j,k)^2 \ell_{kj} \\
&=& w_{\ell_1}^2(-w_{\ell_1})^{-1} + \cdots + w_{\ell_p}^2(-w_{\ell_p})^{-1}\\
&=& -(w_{\ell_1} + \cdots + w_{\ell_p}) \\
&=& - \hat{\delta}_j.
\end{eqnarray*}
Since the $(j,j)$-element of $2D\tilde{\tau} - {\bf 1}\hat{\delta}'$ is $-\hat{\delta}_j,$
the proof is complete in this case.

\vsp
\noindent
{\it Case (ii).} $i \not = j.$
We assume, without loss of generality, that the $ij$-path passes through
$u_1$ (it is possible that $i = u_1$). Let $d(i, j) = \alpha.$
Then $d(i,u_1) = \alpha - w_{\ell_1}, d(i,u_2) = \alpha +  w_{\ell_2},
\ldots, d(i,u_p) = \alpha + w_{\ell_p}.$
 We have
\begin{eqnarray*}
(\Delta L)_{ij} &=& \sum_{k=1}^n d(i,k)^2 \ell_{kj} \\
&=& d(i,u_1)^2(-w_{\ell_1})^{-1} + \cdots + d(i,u_p)^2(-w_{\ell_p})^{-1}+ d(i,j)^2\ell_{jj}\\
&=& (\alpha - w_{\ell_1})^2(-w_{\ell_1})^{-1}  + (\alpha + w_{\ell_2})^2(-w_{\ell_2})^{-1} +\cdots + 
(\alpha + w_{\ell_p})^2(-w_{\ell_p})^{-1}  \\
&+&\alpha^2((w_{\ell_1})^{-1} + \cdots + (w_{\ell_p})^{-1}) \\
&=& (-2\alpha w_{\ell_1} + w_{\ell_1}^2)(-w_{\ell_1})^{-1} + (2\alpha w_{\ell_2} + w_{\ell_2}^2) (-w_{\ell_2})^{-1}
+ \cdots \\
&+& (2\alpha w_{\ell_p} + w_{\ell_p}^2)(-w_{\ell_p})^{-1} \\
&=& 2 \alpha - 2\alpha (p-1) - (w_{\ell_1} + \cdots + w_{\ell_p})\\
&=& 2 \alpha \tau_j -   (w_{\ell_1} + \cdots + w_{\ell_p}),
\end{eqnarray*}
which is the $(i,j)$-element of  $2D\tilde{\tau} - {\bf 1}\hat{\delta}'$
and the proof is complete. \qed

\section{Determinant}

Our next objective is to obtain  a formula for the determinant of the squared
distance matrix.  We first consider the case when the tree has no vertex of 
degree $2.$

\begin{theorem}  \label{lem332}
Let $T$ be a tree with vertex set $V(T) = \{1, \ldots, n\},$
edge set $E(T) = \{e_1, \ldots, e_{n-1}\},$ and edge weights $w_1, \ldots, w_{n-1}.$
Suppose  $T$ has no vertex
of degree $2.$ 
Then
\begin{equation} \label{eq366}
\det \Delta = (-1)^{n-1} \frac{4^{n-2}}{2} \prod_{i=1}^n \tau_i 
\prod_{i=1}^{n-1} w_i^2 \sum_{i=1}^n
\frac{\hat{\delta}_i^2}{\tau_i}.
\end{equation}
\end{theorem}

\noindent  {\bf Proof} We assign an orientation to the edges of the tree and
let  $ H$ and $Q$ be,
respectively, 
edge orientation matrix and the vertex-edge incidence matrix of $T.$

Let $\Delta_i$ denote the $i$-th column of $\Delta,$ and let $t_i$ be the column vector
with $1$ at the $i$-th place and zeros elsewhere, $i = 1, \ldots, n.$
Then
\begin{equation} \label{eq37}
\left [ \begin{array}{c} Q' \\ t_1' \end{array} \right ]
\Delta \left [ \begin{array}{cc}  Q & t_1 \end{array} \right ]
=
\left [ \begin{array}{cc}  Q' \Delta Q  & Q' \Delta_1 \\
\Delta_1' Q & 0  \end{array} \right ].
\end{equation}
Since $\det  \left [ \begin{array}{c} Q' \\ t_1' \end{array} \right ] = \pm 1,$ it follows from
(\ref{eq37}) that
\begin{eqnarray}
\det \Delta &=& \left [ \begin{array}{cc}  Q' \Delta Q  & Q' \Delta_1 \\
\Delta_1' Q & 0  \end{array} \right ] \nonumber \\
&=& \left [ \begin{array}{cc}  -2FHF & Q' \Delta_1 \\
\Delta_1' Q & 0  \end{array} \right ] \mbox{ by Lemma \ref{lem331}} \nonumber \\
&=& (\det (-2FHF))(-\Delta_1'Q(-2FHF)^{-1}Q' \Delta_1) \nonumber \\
&=&  (-2)^{n-1} \prod_{i=1}^{n-1} w_i^2 (\det H)2  \Delta_1'QF^{-1}H^{-1}F^{-1}Q' \Delta_1 \nonumber \\
&=& (-1)^{n-1}2^n \prod_{i=1}^{n-1} w_i^2 (\det H) \Delta_1' QF^{-1}Q' \hat{\tau} Q F^{-1} Q' \Delta_1, \label{eq38}
\end{eqnarray}
in view of Theorem \ref{lem331}.

By Lemma \ref{lem34} we have
\begin{eqnarray}
\Delta_1' QF^{-1}Q' \hat{\tau} Q F^{-1} Q' \Delta_1
&=& \sum_i (2d_{1i}\tau_i - \hat{\delta}_i)^2 \frac{1}{\tau_i} \nonumber \\
& =& \sum_i (4 d_{1i}^2 \tau_i^2 + \hat{\delta}_i^2 - 4 d_{1i} \tau_i \hat{\delta}_i) \frac{1}{\tau_i} \nonumber \\
&=& \sum_i 4 d_{1i}^2 \tau_i + \sum_i \frac{\hat{\delta}_i^2}{\tau_i} - 4 \sum_i d_{1i} \hat{\delta_i} \label{eq39}
\end{eqnarray}

It follows from (\ref{eq39}) and Lemma \ref{lem33} that
\begin{equation} \label{eq40}
\Delta_1' QF^{-1}Q' \hat{\tau} Q F^{-1} Q' \Delta_1 = \sum_i \frac{\hat{\delta}_i^2}{\tau_i}.
\end{equation}
Also by Theorem \ref{lem331},
\begin{equation} \label{eq391}
 \det H = 2^{n-2} \prod_{i=1}^n \tau_i.
\end{equation}
The proof is complete by substituting (\ref{eq40}) and (\ref{eq391}) in (\ref{eq38}). \qed

\begin{corollary}  \cite{krish1} \label{lem333}
Let $T$ be an unweighted  tree with vertex set $V(T) = \{1, \ldots, n\}.$
Suppose  $T$ has no vertex
of degree $2.$ 
Then
\begin{equation} \label{eq36}
\det \Delta = (-1)^{n} {4^{n-2}} \left (2n - 1 - 2 \sum_i \frac{1}{\tau_i}\right )\prod_{i=1}^n \tau_i .
\end{equation}
\end{corollary}

\noindent  {\bf Proof}
We set $w_i = 1, i = 1, \ldots, n-1$ in Theorem \ref{lem332}. Then $\hat{\delta}_i
= \delta_i = 2 - \tau_i, i = 1, \ldots, n.$ We have
\begin{eqnarray}
\sum_i \frac{{\delta_i}^2}{\tau_i} &=& \sum_i \frac{(2-\tau_i)^2}{\tau_i} \nonumber \\
&=& \sum_i \frac{4 + \tau_i^2 - 4 \tau_i}{\tau_i} \nonumber \\
&=& 4 \sum_i \frac{1}{\tau_i} + \sum_i \tau_i - 4n \nonumber \\
&=& 4 \sum_i \frac{1}{\tau_i} + 2 - 4n \nonumber \\
&=& -2 \left (2n - 1 - 2 \sum_i \frac{1}{\tau_i} \right ). \label{eq392}
\end{eqnarray}
The proof is complete by substituting (\ref{eq392}) in (\ref{eq366}). \qed

\vsp

We turn to the case when there is a vertex of degree $2.$

\begin{theorem} \label{lem334}
Let $T$ be a tree with vertex set $V(T) = \{1, \ldots, n\},$
edge set $E(T) = \{e_1, \ldots, e_{n-1}\},$ and edge weights $w_1, \ldots, w_{n-1}.$
Let $q$ be a vertex of degree $2$ and let  $p$ and $r$ be neighbors of $q.$
Let $e_i = (pq), e_j = (qr).$
Then
\begin{equation} \label{eq367}
\det \Delta = (-1)^{n-1} 2^{2n-5} (w_i+w_j)^2  \prod_{s=1}^{n-1} w_s^2
\prod_{k \not = q} \tau_k.
\end{equation}
\end{theorem}

\noindent  {\bf Proof} We assume, without loss of generality, that $e_i$ is directed from $p$
to $q$ and $e_j$ is directed from $q$ to $r.$

\[\xymatrix{
\circ \mbox{$p$} \ar[r]^{e_i} & \circ \mbox{$q$}  \ar[r]^{e_j}  & \circ \mbox{$r$}   
}\]

Let $z_q$ be the $n \times 1$ unit vector with $1$ at the $q$-th place and zeros elsewhere.
Let $\Delta_q$ be the $q$-th column of $\Delta.$  We have
\begin{equation} \label{eq368}
\left [ \begin{array}{c}
Q' \\ z_q' \end{array} \right ] \Delta
\left [ \begin{array}{cc}
Q & z_q \end{array} \right ] =
\left [ \begin{array}{cc}
Q' \Delta Q  & Q' \Delta_q \\
\Delta_q' Q & 0 \end{array} \right ] =
\left [ \begin{array}{cc}
-2FHF  & Q' \Delta_q \\
\Delta_q' Q & 0 \end{array} \right ],
\end{equation}
in view of Lemma \ref{lem34}. It follows from (\ref{eq368}) that
\begin{equation} \label{eq369}
\left [ \begin{array}{cc}
F^{-1} & 0 \\
0 & 1 \end{array} \right ]
\left [ \begin{array}{c}
Q' \\ z_q' \end{array} \right ] \Delta
\left [ \begin{array}{cc}
Q & z_q \end{array} \right ]
\left [ \begin{array}{cc}
F^{-1} & 0 \\
0 & 1 \end{array} \right ]
 =
\left [ \begin{array}{cc}
-2H  & F^{-1}Q' \Delta_q \\
\Delta_q' QF^{-1} & 0 \end{array} \right ].
\end{equation}
Taking determinants of matrices in (\ref{eq369}) we get
\begin{equation} \label{eq370}
(\det F^{-1})^2 \det \Delta = \det \left [ \begin{array}{cc}
-2H  & F^{-1}Q' \Delta_q \\
\Delta_q' QF^{-1} & 0 \end{array} \right ].
\end{equation}
Note that the $i$-th and the $j$-th columns of $H$ are identical. 

Let  $H(j|j)$ denote the submatrix obtained by deleting row $j$ and column $j$ from $H.$
In $\left [ \begin{array}{cc}
-2H  & F^{-1}Q' \Delta_q \\
\Delta_q' QF^{-1} & 0 \end{array} \right ],$ subtract column $i$ from column $j,$
row $i$ from row $j,$ and then expand the determinant along column $j.$ Then we get 
\begin{eqnarray}
 \det \left [ \begin{array}{cc}
-2H  & F^{-1}Q' \Delta_q \\
\Delta_q' QF^{-1} & 0 \end{array} \right ] &=& -((\Delta_q' QF^{-1}))_j - (\Delta_q' QF^{-1})_j)^2 \det (-2H(j|j)) \nonumber \\
&=& -(-2)^{n-2} \det H(j|j) (-w_j - w_i)^2, \label{eq371}
\end{eqnarray}

Note that $H(j|j)$ is the edge orientation matrix of the tree obtained by deleting vertex $q$
and replacing edges $e_i$ and $e_j$
by a single edge directed from $p$ to $r$ in the tree.  Hence by Theorem \ref{lem331},
\begin{equation} \label{eq372}
\det H(j|j) = 2^{n-3} \prod_{k \not = q}\tau_k.
\end{equation}
It follows from (\ref{eq369}),(\ref{eq370}) and (\ref{eq371}) that
\begin{eqnarray}
\det \Delta  &=& -(\det F)^2 (-1)^n 2^{n-2} 2^{n-3} (\prod_{k \not = q}\tau_k)
(w_i+w_j)^2  \nonumber \\
&=&  (-1)^{n-1} 2^{2n-5} (w_i+w_j)^2  \prod_{s=1}^{n-1} w_s^2
\prod_{k \not = q} \tau_k,
\end{eqnarray}
and the proof is complete. \qed 

\vsp

\begin{corollary}  
Let $T$ be a tree with vertex set $V(T) = \{1, \ldots, n\},$
edge set $E(T) = \{e_1, \ldots, e_{n-1}\},$ and edge weights $w_1, \ldots, w_{n-1}.$
Suppose $T$ has at least two vertices of degree $2.$ Then
$\det \Delta = 0.$
\end{corollary}

\noindent  {\bf Proof} The result follows from Theorem \ref{lem334} since $\tau_i = 0$
for at least two values of $i.$ \qed

\section{Inverse}

We now turn to the inverse of $\Delta,$ when it exists. When the tree has no vertex of degree
$2,$ we can give a concise formula for the inverse. We first prove some  preliminary results.

\begin{lemma} \label{lem335}
Let the tree have no vertex of degree $2.$ Then
\begin{equation} \label{eq3721}
\Delta (2\tau - L \hat{\tau} \hat{\delta}) = (\hat{\delta}' \hat{\tau} \hat{\delta}) {\bf 1}.
\end{equation}
\end{lemma}

\noindent  {\bf Proof}
By Lemma \ref{lem35}, $\Delta L = 2 D \tilde{\tau} - {\bf 1} \hat{\delta}'.$
Hence
\begin{equation} \label{eq373}
 \Delta L \hat{\tau} \hat{\delta}  = 2 D \hat{\delta} -(\hat{\delta}' \hat{\tau} \hat{\delta}) {\bf 1}.
\end{equation}
Since by Lemma \ref{lem33}, $\Delta \tau = D \hat{\delta},$ we obtain the result
from (\ref{eq373}). \qed

\vsp
For a square matrix $A,$ we denote by $\cof A,$ the sum of the cofactors of $A.$

\begin{lemma} \label{lem336}
Let $T$ be a tree with vertex set $V(T) = \{1, \ldots, n\},$
edge set $E(T) = \{e_1, \ldots, e_{n-1}\},$ and edge weights $w_1, \ldots, w_{n-1}.$
Suppose $T$ has  no vertex of degree $2.$ Then
\begin{equation} \label{eq374}
{\cof} \Delta  = (-1)^{n-1} 2^{2n-3} \prod_{k=1}^{n-1} w_k^2 \prod_{i=1}^n \tau_i.
\end{equation}
\end{lemma}

\noindent  {\bf Proof} By Lemma \ref{lem34}, $Q' \Delta Q = -2FHF.$ Taking determinant of both sides and using Cauchy-Binet formula,
we get
\begin{eqnarray}
\cof \Delta &=& (-2)^{n-1} (\det F)^2 \det H \nonumber \\
&=& (-2)^{n-1} \prod_{k=1}^{n-1} w_k^2 2^{n-2} \prod _{i=1}^n \tau_i \mbox{ by Theorem } \ref{lem331}\nonumber \\
&=& (-1)^{n-1} 2^{2n-3} \prod_{k=1}^{n-1} w_k^2 \prod_{i=1}^n \tau_i,
\end{eqnarray}
and the proof is complete. \qed

\begin{corollary} \label{lem338}
Let the tree have no vertex of degree $2$ and let $\beta = \hat{\delta}' \hat{\tau} \hat{\delta}.$
If $\beta \not = 0,$ then $\Delta$ is nonsingular and 
\begin{equation} \label{eq375}
{\bf 1}' \Delta^{-1} {\bf 1} = \frac{4}{\beta}.
\end{equation}
\end{corollary}

\noindent  {\bf Proof}
Observe that $\beta = {\displaystyle \sum_{i=1}^n}
\frac{\hat{\delta}_i^2}{\tau_i}.$
By Theorem \ref{lem332},
\begin{equation} \label{eq376}
\det \Delta = (-1)^{n-1} \frac{4^{n-2}}{2} \prod_{i=1}^n \tau_i 
\prod_{i=1}^{n-1} w_i^2  \sum_{i=1}^n
\frac{\hat{\delta}_i^2}{\tau_i}.
\end{equation}
If $\beta \not = 0,$ then $\Delta$ is nonsingular by (\ref{eq376}).
Note that ${\bf 1}' \Delta^{-1} {\bf 1} = \frac{\cof \Delta}{\det \Delta}.$ The
proof is complete using Lemma \ref{lem336} and (\ref{eq376}). \qed

\begin{theorem} \label{lem337}
Let the tree have no vertex of degree $2$ and let $\beta = \hat{\delta}' \hat{\tau} \hat{\delta}.$
Let  $\eta = 2\tau - L \hat{\tau} \hat{\delta}.$ If $\beta \not = 0,$ then $\Delta$ is nonsingular and 
\begin{equation} \label{eq377}
\Delta^{-1} = -\frac{1}{4} L \hat{\tau} L + \frac{1}{4\beta} \eta\eta'.
\end{equation}
\end{theorem}

\noindent  {\bf Proof}
Let $X = -\frac{1}{4} L \hat{\tau} L + \frac{1}{4\beta} \eta\eta'.$ Then
\begin{equation} \label{eq378}
\Delta X = -\frac{1}{4} \Delta L \hat{\tau} L - \frac{1}{4\beta} \Delta \eta\eta'.
\end{equation}

By Lemma \ref{lem35}, $\Delta L = 2 D \tilde{\tau} - {\bf 1} \hat{\delta}'.$
Hence
\begin{equation} \label{eq379}
 \Delta L \hat{\tau} L   = 2 D L - {\bf 1}\hat{\delta}' \hat{\tau}L.
\end{equation}
Using Theorem \ref{lem32}, we can see that
\begin{equation} \label{eq380}
DL = -2I + {\bf 1}\tau'.
\end{equation}
Finally, by Lemma  \ref{lem335}, $\Delta \eta = \beta.$ This fact and (\ref{eq378}),
(\ref{eq379}) and (\ref{eq380})  lead to
\begin{equation} \label{eq381}
\Delta X = I - \frac{1}{2}{\bf 1}\tau' + \frac{1}{4} \hat{\delta}' \hat{\tau} L + \frac{1}{4\beta} {\bf 1}\eta'.
\end{equation}
Since $\eta = 2\tau - L \hat{\tau} \hat{\delta},$ it follows from (\ref{eq381}) that $\Delta X = I$ and
the proof is complete. \qed

\vsp

We conclude with an example to show that the condition $\beta \not = 0$
is necessary in Theorem \ref{lem337}.

\vsp

\noindent
{\bf Example} Consider the tree
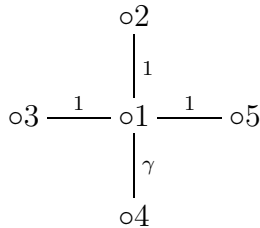
\[\xymatrix{& \circ \mbox{2} \ar@{-}[d]^{1} & & &  & \\
\circ \mbox{3} \ar@{-}[r]^{1}& \circ \mbox{1}\ar@{-}[d]^{\gamma}\ar@{-}[r]^{1} &\circ \mbox{5} &   & \\
&  \circ \mbox{4}& & & 
}\]

The distance matrix of the tree is given by
$$D = \left [ \begin{array}{rrrrr}
0 & 1 & 1 & \gamma & 1\\
1 & 0 & 2 & 1+\gamma  & 2\\
1 & 2 & 0 & 1+\gamma  & 2\\
\gamma  & 1 + \gamma  & 1 + \gamma  & 0 & 1 + \gamma \\
1 & 2 & 2 & 1+ \gamma  &  0\\
\end{array} \right ].$$

It can be checked that $\det \Delta = -32 \gamma^2(\gamma^2 - 6\gamma  - 3).$ Thus $\Delta$ is singular
if $\gamma = 3 + 2\sqrt{3}.$
Note that $\hat{\delta}' = [\gamma+3, 1, 1, \gamma, 1],  \, \tau' = [ -2, 1, 1, 1, 1]$
and hence, if  $\gamma = 3 + 2\sqrt{3},$ then $\displaystyle{\sum_{i=1}^4} \frac{\hat{\del}^2}{\tau_i} = 0.$

\vsp
{\bf Acknowledgment} I sincerely thank Ranveer Singh for a careful reading of the manuscript.
Support from the JC Bose Fellowship, Department of Science and Technology, Government
of India, is gratefully acknowledged.

\vsp

\end{document}